%% file: paper3.tex
\documentclass[12pt,twoside]{gen-j-l}
\usepackage{amsthm}
\usepackage{amsmath}
\usepackage{amstext}
\usepackage{amsfonts}
\usepackage{latexsym}
\usepackage{amscd}
\usepackage{diagrams}
\newcommand{\la}{\ensuremath{\longrightarrow}}

\newcommand{\p}{\ensuremath{\mathbb{P}}}
\newcommand{\pone}{\ensuremath{\mathbb{P}^{1}}}

\newcommand{\I}{\ensuremath{\mathcal{I}}}

\newcommand{\sheaf}{\ensuremath{\mathcal{O}}}

\newtheorem{theorem}{Theorem}[section]
\newtheorem{lema}[theorem]{Lemma}
\newtheorem{definition}[theorem]{Definition}
\newtheorem{proposition}[theorem]{Proposition}
\newtheorem{corollary}[theorem]{Corollary}

\newtheorem{example}[theorem]{Example}
\newarrow{BirTo}{}{dash}{}{dash}{>}

\begin{document}
\title[3-FOLD CONTRACTIONS, SURFACE TO CURVE I]{TERMINAL 3-FOLD DIVISORIAL CONTRACTIONS 
OF A SURFACE TO A CURVE I}
\author{Nikolaos Tziolas}
\address{Mathematics Institute, University of Warwick, Coventry, CV4 7AL, England}
\email{tziolas@maths.warwick.ac.uk}
\thanks{The author is supported by the Marie Curie research fellowship grant number MCFI-1999-00448.}

%    General info
\subjclass{Primary 14C05, 14J32; Secondary 14D15}

\date{September 2001}

%\dedicatory{This paper is dedicated to our authors.}

\keywords{Algebraic geometry}

\begin{abstract}
Let $\Gamma \subset X$ be a smooth curve on a 3-fold which has only index 1 terminal singularities along this curve. 
In this paper we investigate the existence of extremal terminal 
divisorial contractions $E\subset Y \la \Gamma \subset X$, contracting 
an irreducible surface $E$ to $\Gamma$. We consider cases with respect to 
the singularities of the general hypersurface section $S$ of $X$ through $\Gamma$. 
We completely classify the cases when $S$ is $A_i$, $i \leq 3$, and $D_{2n}$ for 
any $n$.
\end{abstract}

\maketitle

\input{introduction}

\input{section1}

\input{section2}

\input{section3}
\input{section4}

\input{bibliography}
\end{document}

%% file: introduction.tex
\setcounter{section}{-1}
\section{Introduction}
One of the main objectives of birational geometry is to identify in each birational class of varieties some distinguished 
members which are ``simple'' and are called minimal models, and then study the structure of birational maps between them. 

In dimension two, satisfactory results were known for over one hundred years. In higher dimensions, a program called the minimal 
model program (MMP) was developed to search for minimal models. After contributions of Reid, Mori, Kawamata, Koll\'{a}r, Shokurov 
and others, the program was completed in dimension three by Mori in 1988. According to this program, a variety $X$ is called a 
minimal model iff it is $\mathbb{Q}$-factorial, terminal and $K_X$ is nef. According to Mori's theorem:

\begin{theorem}[Mori88]
Let $X$ be a $\mathbb{Q}$-factorial, terminal projective 3-fold. Then there is a sequence of birational maps \[
X=X_1 ---> X_2 --->\cdots --->X_n=X^{\prime}\]
such that $X^{\prime}$ is $\mathbb{Q}$-factorial, terminal projective and exactly one of the following cases happen:
\begin{enumerate}
\item $K_{X^{\prime}}$ is nef and hence $X^{\prime}$ is a minimal model.
\item $X^{\prime}$ is a Mori fiber space. This means that there exists a morphism $f: X^{\prime} \la S$ such that $S$ is normal, 
$\dim S \leq 2$, $\rho (X^{\prime}/S)=1$ and $-K_{X^{\prime}}$ is $f$-ample.
\end{enumerate}
\end{theorem}
The birational maps that appear in the theorem are divisorial contractions and flips. It is also known~\cite{Ko-Mo98} that any 
birational map between minimal models is an isomorphism in codimension one and a composition of flops. Terminal flops were classified 
by the work of Koll\'{a}r~\cite{Ko91}. 

The structure of birational maps between Fano fiber spaces is much more complicated. A program called the Sarkisov program 
was developed by Sarkisov, Reid and Corti to factorize birational maps between these spaces as a composition of so called 
``elementary links''~\cite{Cor95}. These links consist of flops, flips and divisorial contractions. Therefore to better understand 
the structure of birational maps between Fano fiber spaces, it is important to understand divisorial contractions and flips. 
Flips were classified by Koll\'{a}r and Mori~\cite{Ko-Mo92}. However the structure of divisorial contractions is not yet completely 
understood. 

Let $E \subset Y \stackrel{f}{\la} \Gamma \subset X$ be a divisorial contraction. Mori and Cutkoski completely classify such 
contractions when $Y$ is Gorenstein. In particular, if $\dim \Gamma =1$, then $X$ is smooth along $\Gamma$ and $Y$ is just 
the blow up of $X$ along $\Gamma$. 

Moreover, Kawamata~\cite{Kaw94} showed that if there is a point $P \in \Gamma \subset X$ such that $P \in X$ is a cyclic quotient 
terminal singularity, then $\Gamma = \{P\}$ and $f$ is a weighted blow up.

Divisorial contractions of a surface to a point, i.e. when $\Gamma =\{pt\}$, have been studied by Luo, Corti, Kawakita and others. 

This paper studies divisorial contractions of a surface to a curve, i.e. when $\dim \Gamma =1$ and 
$X$ has only index 1 terminal singularities along $\Gamma$. It is not always true that given $\Gamma \subset X$, then there is a terminal 
contraction of a surface to $\Gamma$. We investigate when there is one, give criteria for existence or not and in the case that 
there is a terminal contraction we also describe the singularities of $Y$. 

By Reid's general elephant conjecture~\cite{Ko-Mo92}, if a terminal contraction exists then there is a DuVal section 
$\Gamma \subset S \subset X$. We base the classification of contractions on the type of singularities of a general $S$ through 
$\Gamma$ instead of $X$ itself.

Theorem~\ref{canonicalcase} shows that under certain conditions there is always a canonical contraction. 

The objective of the rest of this paper is to investigate when the contraction is terminal in the case that $\Gamma$ is smooth. 

I would now like to point out that for most applications this is the most interesting case and therefore it is not a big 
restriction to concentrate on the case that $\Gamma$ is smooth. In particular, when one studies the birational rigidity 
of a Fano 3-fold, it is important to exclude certain curves as maximal centers. In most cases~\cite{Cor-Rei00} these curves 
are either lines or conics. Moreover, Miles Reid believes that a curve that is a maximal center can have at worst nodes as 
singularities. This is supported by a large class of examples that he calculated. Finally I believe that the methods 
of this paper can be used to work the case when $\Gamma$ has at worst local complete intersection singularities. 

In order to investigate the existence of a terminal contraction, it is important to obtain normal forms for the equations 
of \mbox{$\Gamma \subset S  \subset X$}. This is done in Proposition~\ref{normalforms}.

Theorem~\ref{An} gives criteria for existence in the case that $S$ is $A_i$ with $i \leq 3$ and 
theorem~\ref{Dn} treats the case when $S$ has $D_{2n}$ singularities. 

There is an important difference between the $D_{2n}$ and $D_{2n+1}$, as well as for the higher $A_n$ cases. The main 
difficulty is the explicit calculation of the $\mathbb{Q}$-factorialization $Z$ of $E_1$, as appears in the proof of 
theorem~\ref{canonicalcase}. The reason of this difficulty becomes clear in lemma 6.2. 
I have developed a method of calculating $Z$ and under one technical restriction 
i was succesfull in obtaining results for the $D_{2n+1}$ cases too. Unfortunately at this moment I do not see how to remove this 
restriction or how to present the results in a ``nice'' and not too technical way. I plan to address these problems in a follow up paper.

Part of this work was done during my visit to Princeton University in April 2001. 

I would like to thank J\'{a}nos Koll\'{a}r for many interesting and useful discussions during my visit to Princeton, 
and Miles Reid who suggested this problem to me.

%% file: section1.tex
\section{Uniqueness and Canonical Contractions}
First we define divisorial contractions.
\begin{definition}
A 3-fold divisorial contraction is a morphism\[
f: E \subset Y \la \Gamma \subset X \]
such that:
\begin{enumerate}
\item $X$ and $Y$ are $\mathbb{Q}$-factorial.
\item $Y-E \cong X-\Gamma$ and $E$ is a prime divisor.
\item $-K_Y$ is $f$-ample. 
\item $\rho (Y/X)=1$.
\end{enumerate}
\end{definition}
The next proposition shows that under certain conditions, the contraction if it exists is unique.
\begin{proposition}
Let $f: E \subset Y \la \Gamma \subset X$ be a 3-fold divisorial contraction of an irreducible 
surface $E$ to a curve $\Gamma$. Suppose that $X$, $Y$ are normal, $\dim f(Y^{sing})=0$, $X$ has 
isolated singularities and $-E$ is $f-$ample. Then 
\[
Y \cong Proj \oplus _{d \geq 0} \I_{\Gamma, X}^{(d)}
\]
\end{proposition}
\begin{proof}
Since $-E$ is $f-ample$, it follows that \[
Y \cong \oplus_{d \geq 0} f_{\ast} \sheaf_{Y}(-dE). 
\]
The proposition will follow once it is shown that $f_{\ast} \sheaf_{Y}(-dE)=\I_{\Gamma, X}^{(d)}$. Since 
by assumption $X$ has isolated singularities and $\dim f(Y^{sing})=0$, there are finitely many points 
$P_1,\ldots,P_k$ in $X$ so that $X-\{P_1,\ldots,P_k\}$ and $Y-f^{-1}\{P^1,\ldots, P_k\}$ are smooth. Hence by 
Mori's classification of smooth terminal contractions, $f$ is just the blow up of $\Gamma$ away from 
$P_1,\ldots,P_k$. Hence \[
f_{\ast} \sheaf_{Y}(-dE) \mid_{X-\{P_1,\ldots,P_k\}} = \I_{\Gamma, X}^{d}\mid_{X-\{P_1,\ldots,P_k\}} 
\]
Now let $ s\in \sheaf_{X}$. Since $\dim f^{-1}(x) \leq 1$, $\forall x \in X$, \[
s\in f_{\ast} \sheaf_{Y}(-dE) \Longleftrightarrow s\mid_{X-\{P_1,\ldots,P_k\}} \in 
\I_{\Gamma, X}^{d}\mid_{X-\{P_1,\ldots,P_k\}}.\]
The proposition will now follow immediately from the following simple property of symbolic powers. 

\begin{lema}
Let $A$ be a ring and $\mathrm{m}_1, \ldots , \mathrm{m}_k$ be maximal ideals. Let $I \subset A$ 
be a prime ideal such that $I \neq \cap_{j=1}^{k} \mathrm{m}_{j}$. Let $a \in A$ such that \[
a \mid_{Spec A - \{\mathrm{m}_1, \ldots , \mathrm{m}_k\} } \in I^d.\]
Then $a \in I^{(d)}$.
\end{lema}
\begin{proof}
$SpecA-\{\mathrm{m}_1, \ldots , \mathrm{m}_k\}=\cup_{s\in \cap \mathrm{m}_{i}} Spec A_{s}$. The 
assumption implies that $a \in I^{d}A_{s}$, $\forall s\in \cap_i \mathrm{m}_{i}$. Let 
$s \in \cap_i \mathrm{m}_{i} - I$. Then $s^{\nu}a \in I^{d}$, for some $\nu$. Hence $a \in I^{(d)}$.
\end{proof}
\end{proof}
In particular, terminal contractions of a surface to a curve are unique. However this is not true for 
canonical if the condition $\dim f(Y^{sing})=0$ is removed, as shown by the next example:
\begin{example}
Let $X$ be given by $x^2 + y^2 z + z^3 + t^5=0$ and $\Gamma : x=z=t=0$. Then there are 2 non-isomorphic 
canonical divisorial contractions $g_i: S_i \subset Z_i \la \Gamma \subset X$ contracting the surfaces 
$S_i \subset Z_i$ to $\Gamma$. $Z_1$ has index 1 and is singular along a section of $S_1 \la \Gamma$ 
and $Z_2$ has index 2 and its singular locus is $g_2^{-1}(0)$. 
\end{example}
\begin{proof}
The statement about the contraction $Z_2 \la X$ follows from Theorem~\ref{Dn}. So here we will only show 
how to construct the contraction \mbox{$g_1 :Z_1 \la X$.} 

Let $I=(x^2, z, t)$. Let $ g_1 : Z_1 = B_I X \la X$ be the blow up of the ideal $I$ in $X$. 
$Z_1 \subset \mathbb{C}^4 \times \mathbb{P}^2$. Let $u$, $v$, $w$ be coordinates for $\mathbb{P}^2$. 
Look at the chart $w \neq 0$. $Z_1$ is given by 
\begin{gather*}
x^2 -ut=0 \\
u+y^2v+v^3t^2+t^4=0
\end{gather*}
There is only one $g_1$-exceptional divisor $S_1$ given by $x=t=u+y^2 u=0$. $Z_1$ is easily seen to be singular 
along the line $l : x=t=u=v=0$ which lies over $\Gamma$. Moreover, since $Z_1$ is a complete intersection, 
it has index 1.
\end{proof}
The singularities of the general hyperplane section $S$ of $X$ through $\Gamma$ are very important for the 
study of terminal contractions as shown by the following theorem.
\begin{theorem}[~\cite{Ko-Mo92}]
Let $f : X \supset C \la Y \ni Q$ be an extremal neighborhood. Then the general member $E_{X}$ 
of $|-K_{X}|$ and $E_{Y}=f(E_{X}) \in |-K_{Y}|$ have only DuVal singularities. Moreover, the minimal resolution 
of $E_{Y}$ dominates $E_{X}$.
\end{theorem}
From this it follows that if a terminal contraction exists, then there is a DuVal section $S$ containing $\Gamma$. 
The converse is not true as shown by the next example. 

\begin{example}
Let $X$ be given by $x^2 + y^2 z + z^n + t ^m =0$, with $n,m \geq 4$, and \mbox{$\Gamma: x=z=t=0$.} Then the contraction exists 
but it is only canonical, even though the section $S$ given by $t=0$ is $D_{n+1}$. 
\end{example}
\begin{proof}
The section $S$ through $\Gamma$ given by $t=0$ is $D_{n+1}$, $\Gamma \subset S$ is of type $DF_l$ (look at definition~\ref{DF}) 
and therefore by Lemma 7.3 so is the general section. Moreover, by lemma 7.4, there is no $D_4$ section of $X$ containing $\Gamma$. 
Now the result follows from Theorem~\ref{Dn}.2.a. 
\end{proof}

The next example shows that there are cases when there is no DuVal section containing $\Gamma$, and hence we can 
immediately conclude that there is no terminal contraction.

\begin{example}
Let $X$ be given by $x^2 + y^3 + z^3 + y t^6 =0$ and $\Gamma: x=y=z=0$. Then there is no DuVal section of $X$ containing 
$\Gamma$ and hence there is no terminal contraction. Moreover, the blow up $Y=B_{\Gamma}X$ of $X$ along $\Gamma$ is 
not even canonical. 
\end{example}
\begin{proof}
First observe that $0 \in X$ is a $cD_4$ singularity and therefore terminal. The section $t=0$ is the surface\[
x^2+y^3+z^3=0\]
which is easily seen to have $D_4$ singularities.

A general hyperplane containing $\Gamma$ in $\mathbb{C}^4$ is given by $x=ay+bz$. Hence the corresponding section 
of $X$, $S$, is \[
S \; : \; (ay+bz)^2+y^3+z^3+yt^6=0.\]
This is not DuVal. To see this blow up the origin. In the affine chart given by $y=yt$, $z=zt$, the blow up is given by  \[
(ay+bz)^2 +y^3t+z^3t+yt^5=0.\]
Blow up again. In the affine chart $y=yt$, $z=zt$, the blow up is given by \[
(ay+bz)^2+y^3t^2+z^3t^2+yt^4=0.\]
This is singular along the line $ay+bz=t=0$ and hence it is not normal. Therefore the original singularity is not DuVal. Hence there 
is no DuVal section of $X$ containing $\Gamma$. 

Now let $f : Y=B_{\Gamma}X \la X$ be the blow up of $X$ along $\Gamma$. In the affine chart $x=xz$, $y=yz$, $Y$ is given by\[
Y \; : \; x^2 z +y^3 z^2 +z^2 + yt^6=0.\]
Moreover, \[
x^2 z +y^3 z^2 +z^2 + yt^6=x^2 z +(y^3 +1)z^2 + yt^6 \sim x^2z+z^2+yt^6 \]
which is singular along the line $l: x=z=t=0$. A typical section given by $y=a$ is , \[
x^2z+z^2+at^6=0. \]
Blow up the origin twice as before to get a nonnormal surface. Hence $Y$ has a line of at best log canonical singularities. 

Hence the existence of a DuVal section through $\Gamma$ is important to conclude that $Y=B_{\Gamma} X$ is canonical in 
Proposition~\ref{blowup}.
\end{proof}
The previous example shows the significance of the existence of a DuVal section $S$ of $X$ through $\Gamma$. 
Therefore from now on we will assume the existence of such a section. In fact, we will study the existence of 
terminal contractions by considering cases with respect to the type of singularities of the general $S$ through $\Gamma$, 
instead of the singularities of $X$ itself.

To start with we will show that there is always a canonical contraction. 

\begin{theorem}\label{canonicalcase}
Let $X$ be a $\mathbb{Q}$-factorial index 1 terminal 3-fold and $\Gamma \subset X$ an irreducible curve 
having at worst plane curve singularities. Suppose there is a DuVal section $S$ of $X$ containing $\Gamma$. 
Then there is an extremal contraction $g : W \la X$ 
contracting an irreducible surface $E$ to $\Gamma$, with $W$ canonical and $\dim g(W^{sing})=0$. In particular, 
$R(\Gamma,X)=\oplus_{n\geq 0}\I^{(n)}_{\Gamma,X}$ is finitely generated.
\end{theorem}
The proof of the previous theorem will be given in section 3.

\section{Some easy lemmas}
\begin{definition}
Let $X$ be a normal algebraic variety, and $D$ a $\mathbb{Q}$-Cartier divisor in $X$. 
\begin{enumerate}
\item Let $P \in X$. Then the index of $D$ at $P$, $\text{index}_P(D)$, is defined to be the smallest 
$r \in \mathbb{N}-{0}$ such that $rD$ is Cartier at $P$. 
\item The global index of $D$ in $X$, $\text{index}_X(D)$ is the smallest $r \in \mathbb{N}-{0}$ such that $rD$ is Cartier.
\end{enumerate}
\end{definition}
\begin{lema}
Let $X$ be a threefold. Suppose that its singular locus is an irreducible curve $\Gamma$ and that it has only 
hypersurface singularities. Let $D$ be a $\mathbb{Q}$-Cartier divisor on $X$. Then \[
\text{index}_{P}(D)=\text{index}_{X}(D) \]
for any point $P \in \Gamma$. Hence the index of $D$ can be computed at any point of $\Gamma$. 
\end{lema}
\begin{proof}
The proof is almost identical to that for the case of isolated index 1 terminal singularities that appears in~\cite{Kaw88}. I am 
not aware of a reference for this more general case and therefore I include it for the convenience of the reader.

Let $r=index_{X}(D)$. Then there are finitely many possibilities for $index_{P}(D)$. Suppose that $r_1 \leq r_2 \leq \cdots \leq r_k=r$ 
be these possibilities. Suppose that $r_1=index_{P_1}(D)$. Since $\Gamma$ is irreducible, $r_1 D$ is Cartier at all but finitely 
many points where $D$ has index greater than $r_1$. We can assume that it is only one, say $P$, since the result is local. Hence 
$index_{P}(D)=r$. Then $D^{\prime}=r_1 D$ is Cartier everywhere except $P$. Let \[
\pi : W \la X \]
be the index 1 cover of $D^{\prime}$. Hence \[
W-\pi ^{-1}(P) \la X-P \]
is \'{e}tale. $X$ has hypersurface singularities and therefore by~\cite[Th. 5.2]{Mil68} \[
\pi _1 (X-P) = 0.\]
Hence the cover is trivial and therefore $D^{\prime}$ is Cartier.
\end{proof}
\begin{lema}
Let $P \in \Gamma \subset X$. Suppose that $P \in X$ is a 3-dim index 1 terminal singularity, and that $P \in \Gamma$ is 
a plane curve singularity. Let \mbox{$ f : Y=B_{\Gamma}X \la X$} be the blow up of $X$ along $\Gamma$. 
Then $f^{-1}(\Gamma)=E_1 + m E_2$. $E_1$ is an irreducible surface over $\Gamma$. If $P \in \Gamma$ is smooth, then 
$E_2 \cong \mathbb{P}^2$. Otherwise $E_2 \cong \mathbb{P}^2$ or $E_2=\emptyset$.
\end{lema}
\begin{proof}
The result is local around $P$. Since $P \in X$ is an index 1 terminal singularity, $P$ is a cDV point. 
Hence we can assume that $X \subset \mathbb{C}^4$. Suppose $P \in \Gamma$ 
is smooth. First we will show that $\dim f^{-1}(P)=2$. Suppose not. Let $P \in S \subset X$ be a general 
hypersurface section transversal to $\Gamma$. Then $P \in S$ is DuVal and $S^{\prime}=f_{\ast}^{-1}S$ is just 
the blow up of $P$ in $S$. In particular, it is normal and $K_{S^{\prime}}=f^{\ast}K_{S}$. If $\dim f^{-1}(P) =1 $, 
then $f^{\ast}S = S^{\prime}$. In particular, $S^{\prime}$ is Cartier. Moreover, since $f$ is generically 
the blow up of a smooth curve, $K_{Y}=f^{\ast}K_{X}+E_1$. Since $Y$ is CM and $S^{\prime}$ Cartier, 
\begin{gather*}
K_{S^{\prime}}=(K_{Y}+S^{\prime})\mid _{S^{\prime}} = (f^{\ast}K_{X} + E_1 + f^{\ast}S)\mid_{S^{\prime}}= \\
= (f^{\ast}(K_{X}+S)+ E_1)\mid_{S^{\prime}}=f^{\ast}K_{S} + E_1\mid_{S^{\prime}}=K_{S^{\prime}}+ E_1\mid_{S^{\prime}}
\end{gather*}
It now follows that $E_1\mid_{S^{\prime}}=0$, which is impossible. Hence $\dim f^{-1}(P)=2$. Now in both cases,
since $\Gamma$ has at worst plane curve singularites at $P$, $\I_{\Gamma, \mathbb{C}^4}$ is generated by a regular sequence 
$\{g_1 , g_2 , g_3\}$. Hence $B_{\Gamma} \mathbb{C}^4 \la \mathbb{C}^4$ is the blow up of a regular sequence and 
thus all its fibers over $\Gamma$ are isomorphic to $\mathbb{P}^2$. Hence $E_2 \cong \mathbb{P}^2$.
\end{proof}
\begin{lema}
Let $P\in \Gamma \subset S$. Suppose that $P \in \Gamma$ is at worst a plain curve singularity, 
$\Gamma -P$ is smooth, and $S$ is a normal and canonical surface. 
Let $f : S^{\prime}=B_{\Gamma}S \la S$ be the blow up of $S$ along $\Gamma$. Then $S^{\prime}$ is normal.
\end{lema}
\begin{proof}
Clearly, $S^{\prime}-f^{-1}(P) \cong S-P$ since it is the blow up of a Cartier divisor. So the result is local over $P$. 
Hence we can assume that \mbox{$S \subset X=\mathbb{C}^3$}. Let $f : Y=B_{\Gamma}X \la X $ be the blow up of $X$ along $\Gamma$, 
and $E$ the $f-$exceptional divisor. Since $P \in \Gamma$ is a lci singularity, $Y$ is just the blow up of a regular 
sequence, say $\{g_1 , g_2 \}$ in $X$. Hence \[
E \cong \mathbb{P}_{\Gamma}(I_{\Gamma,X} / I_{\Gamma,X}^2).
\]
In particular, $\dim f^{-1}(x) \leq 1 $ $\forall x \in X$. Hence since $Y$ is the blow up of a smooth curve away from $P$, 
$Y$ is normal. Moreover, \[
\begin{array}{lcr}
f^{\ast}S=S^{\prime} + E    & \text{and} &  K_{Y}=f^{\ast}K_{X} +E 
\end{array}.\]
Hence \[
K_{Y}+S^{\prime}=f^{\ast}(K_{X}+S).\]
By~\cite[Th. 7.3]{Ko97}, the pair $(X,S)$ is also canonical. Hence $(Y,S^{\prime})$ is also canonical and hence $plt$. 
By~\cite[Th. 5.51]{Ko-Mo98} and~\cite{Ko97}, $S^{\prime}$ is normal and canonical.
\end{proof}

\begin{lema}
Let $P \in \Gamma \subset X$. Assume that $\Gamma$ is smooth and that $P \in X$ is a 3-dimensional normal hypersurface 
singularity. Let $ f : Y=B_{\Gamma}X \la X$ be the blow up of $X$ along $\Gamma$. Let $P \in S \subset X$ be a general 
hypersurface section through $P$. Then \[
f^{\ast}S=S^{\prime}+(m_{P}X - 1)E_2 \]
where $ m_{P}S$ is the multiplicity of $S$ at $P$ and $f^{-1}(\Gamma)=E_{1}+dE_{2}$ as in lemma 2.3. 
\end{lema}
\begin{proof}
Suppose that $f^{\ast}S=S^{\prime}+aE_2$. $f^{-1}(\Gamma)=E_{1}+dE_{2}$ is Cartier and $\sheaf_{Y}(E_{1}+dE_{2})=\sheaf_{Y}(-1)$. 
Then, 
\[
f^{\ast}S \cdot (E_1+dE_2)^2 = S^{\prime} \cdot (E_1+dE_2)^2 + a E_2 \cdot (E_1+dE_2)^2 \]
Since $E_2 \cong \mathbb{P}^2$, it follows that $E_2 \cdot (E_1+dE_2)^2=1$. Moreover, $S^{\prime}=B_{P}S$ the blow up 
of $S$ at $P$. Let $F$ be the exceptional divisor. Then $F^2 = -m_{P}S = m_{P}X$, and $F^2=S^{\prime}\cdot (E_1+dE_2)^2$. 
Moreover, $(E_1+dE_2)^2=E_1 \cdot (E_1+dE_2)+dE_2 \cdot (E_1+dE_2)= \Delta + dL$, where $\Delta$ is a section of 
$E_1 \la \Gamma$ and $L$ a line in $E_2=\mathbb{P}^2$. Hence\[
f^{\ast}S \cdot (E_1+dE_2)^2 =S \cdot \Gamma =1 \]
and thus \[
a=m_{P}S-1=m_{P}X-1 \]
\end{proof}
\begin{lema}
Let $f : E \la \Gamma$ be a morphism from an irreducible surface to an irreducible curve. Suppose that $f^{-1}(P_0)=C_0$ 
is an irreducible curve and $P_0 \in \Gamma$ a smooth point. Let $f^{-1}(P)=\sum_i a_i C_i $ for another point. 
Then $C_0 \equiv \sum_i b_i C_i$ for $b_i \geq 0$.
\end{lema}
\begin{proof}
Let $\Gamma ^{\prime} \la \Gamma$ be the normalization of $\Gamma$, and 
$E^{\prime}= \Gamma ^{\prime} \times_{\Gamma} E$. Hence there is a fiber square \[
\begin{CD}
E^{\prime}  @>{p^{\prime}}>> E \\
@V{g}VV            @VV{f}V \\
\Gamma ^{\prime} @>{p}>> \Gamma
\end{CD} \]
Since $\Gamma ^{\prime}$ is smooth, for any point $Q$ over $P$, $g^{-1}(P) \equiv g^{-1}(Q)$. Apply $p^{\prime}_{\ast}$ 
to get the result.
\end{proof}
\begin{proposition}\label{blowup}
Let $P \in \Gamma \subset X$. Assume that $P\in X$ is an index 1 terminal singularity and $X$ normal. 
let $ f : Y=B_{\Gamma}X \la X$. Then:
\begin{enumerate}
\item If $\Gamma$ is smooth, then $Y$ is normal of index 1.
\item If there is a DuVal section $P \in \Gamma \subset S \subset X$ and $\Gamma$ has at worst plane curve singularities, 
then $Y$ is normal and canonical of index 1.
\item $K_{Y}=f^{\ast}K_{X}+E_1 + d E_2$, and $f^{-1}(\Gamma)=E_1 + dE_2$ is Cartier but $E_1$, $E_2$ are not 
$\mathbb{Q}$-Cartier.
\end{enumerate}
\end{proposition} 
\begin{proof}
First we will show that $Y$ is CM and that $\omega_Y$ is invertible. The result is local around $P$. 
Since $P \in X$ is index 1 terminal, it is cDV. So we can assume that $X \subset \mathbb{C}^4=U$. Let $W=B_{\Gamma}\mathbb{C}^4$. 
Since $\Gamma$ has at worst plane curve singularities, $W$ is the blow up of a regular sequence and 
therefore $W$ is lci and hence CM and $\omega_W$ is invertible. 
Moreover, 
\[
E=f^{-1}(\Gamma)=E \cong \mathbb{P}_{\Gamma}(I_{\Gamma,U} / I_{\Gamma,U}^2)
\] 
is Cartier and irreducible. Suppose that $f^{\ast}X=Y+aE$, $a \in \mathbb{N}$. Hence $Y \subset W$ is Cartier and hence is CM and 
$\omega_Y$ is invertible.

Now assume that $\Gamma$ is smooth. Since $P \in X$ is cDV, $m_P X=2$. Let 
$P\in S \subset X$ be a general hypersurface section. Then $P \in S$ is DuVal. From lemma 2.5 it follows that 
\[
f^{\ast}S=S^{\prime}+E_2 
\]
Since $E_2 \cong \mathbb{P}^2$ and $ S^{\prime}+E_2$ is Cartier, it follows that $Y$ is smooth at the generic point of 
$E_2$ and hence regular in codimension 1 and therefore normal. This shows $1$.
 
Now suppose that $f^{-1}(\Gamma)=E_1 + d E_2$. This is of course Cartier. Since $Y \subset W$ 
is Cartier, we can use adjunction to calculate $K_Y$. Since $W-f^{-1}(P)$ is the blow up of a smooth curve, it follows that 
\[
K_W=f^{\ast}K_U + 2E. \]
moreover, $f^{\ast}X=Y+E$. Hence 
\begin{gather*}
K_Y = (K_W + Y)\mid _Y = f^{\ast}(K_U + X)+E \mid_Y = f^{\ast}K_X +E\mid_Y= \\
=f^{\ast}K_X+E_1 + dE_2. 
\end{gather*}
If $E_i$ were $\mathbb{Q}$-Cartier, then for a fiber $\delta$ of $f$ disjoint from $E_2$ and a line $l \subset E_2$, 
$\delta = a l$, $a \geq 0$. Hence $0 \leq l \cdot E_1 = a \delta \cdot E_1 = a \delta \cdot (E_1 + d E_2)= -a <0$, 
which is not possible. Of course one could argue that in this case the exceptional would have to be irreducible 
which is not the case.

Now suppose that a DuVal section $\Gamma \subset S \subset X$ exists. Then $S = X \cap H$ for some general plane in $U$. 
Therefore, $f^{\ast}H=H^{\prime}+E$ and hence it follows that \[
f^{\ast}S=S^{\prime}+E_1 + dE_2. \]
In particular, $S^{\prime}$ is Cartier. From lemma 2.4 it follows that $S^{\prime}$ is normal and canonical and therefore since 
$S^{\prime}$ is Cartier, $Y$ is smooth at some points of $E_2$ and hence regular in codimension 1 and therefore normal.

Since $X$ is terminal and $S$ is canonical it follows that $(X,S)$ is canonical. 
Adjunction for $S^{\prime}$ in $Y$ gives that \[
K_Y + S^{\prime}= f^{\ast}(K_X + S). \]
Hence the pair $(Y,S^{\prime})$ is also canonical. Since $Y$ has index 1 and $S^{\prime}$ is DuVal, it follows that 
$Y$ is also canonical.
\end{proof}
If there is no DuVal section $S$ containing $\Gamma$, then $Y$ may fail to be canonical. In fact as shown by Example 1.7, 
it may not even be log canonical. 

\section{Proof of Theorem~\ref{canonicalcase}}

Let $ f : Y=B_{\Gamma}X \la X$. Then by proposition~\ref{blowup}, $Y$ is canonical and normal. 
Therefore by~\cite{Kaw88}, there is \mbox{$g : Z \la Y$} such that $g$ is an isomorphism in codimension 1 
and $-E_{1}^{Z}=g_{\ast}^{-1}E_1$ is $\mathbb{Q}-$Cartier and $g$-ample. 
More precisely, $R(E_1,Y)= \oplus_{i}\sheaf_{Y}(-iE_1)$ is finitely generated and \mbox{$Z=\text{Proj} R(E_1,Y)$}. 
Hence $K_Z=g^{\ast}K_Y$ and $Z$ is canonical. In fact since $X$ is $\mathbb{Q}$-factorial, $Z$ 
is also $\mathbb{Q}-$factorial. We want to contract 
$E_{2}^{Z}=g_{\ast}^{-1}E_2$ over $X$ and obtain the required contraction. Let $C_i$ be the $g$-exceptional curves. 
Of course $C_i \subset E_1^{Z}$. Then 
$K_Z \cdot C_i =0$. On the other hand, for a $l \subset E_2^{Z}$, the birational transform of a line in $E_2=\mathbb{P}^2$, 
\[
K_Z \cdot l = K_Y \cdot l = (E_1 + dE_2) \cdot l = -1 <0.
\] 
Therefore, $\overline{NE}(Z/X)$ contains $K_Z$-negative extremal rays. Let $R$ be one. Then $R=\mathbb{R}_{+}[\Delta]$ 
for an irreducible curve $\Delta$. $\Delta$ must be contained in $E_2^{Z}$. Suppose not. Then it is a fiber of 
$E_1^{Z} \la \Gamma$. Let $L = E_1 \cap E_2$. Then from lemma 2.6, $\Delta \equiv L^{Z} + \sum_i a_i C_i$, with 
$a_i > 0$. Hence in this case $R$ cannot be extremal. If it was then $[C_i] \in R$ which is impossible. 
So all the $K_Z$-negative extremal curves are contained in $E_2^{Z}$. 
Now by the relative cone theorem there is a contraction $p : Z \la W $ over $X$. There is a factorization:
\[
\begin{CD}
Z @>{g}>> Y \\
@V{p}VV   @VV{f}V \\
W @>{q}>>  X
\end{CD}
\]
By~\cite[Prop. 3.38]{Ko-Mo98}, $W$ is also canonical and for any exceptional divisor, $F$, with center in 
$p(E_2^{Z})$, $a(F,W)>0$. Moreover by~\cite{Benv85}, since $Z$ is Gorenstein, $p : Z \la W $ is divisorial and therefore 
\[
q : W \la X\] 
is the required contraction.

\qed

From the previous proof it becomes clear that $W$ is terminal iff $Z$ has only isolated terminal singularities along the 
$g$-exceptional curves $C_i$ and away from $E_2^{Z}$. To investigate when this happens, 
it is important to obtain an explicit description of the 
$\mathbb{Q}-$factorialization of $E_1$, $Z$. In general this is difficult. However in the cases that the general section 
$S$ of $X$ through $\Gamma$ is $D_{2n}$ or $A_i$ with $i \leq 3$, it is possible to get such a description and therefore treat 
these cases completely.

%% file: section2.tex
\section{Normal forms for $\Gamma \subset S \subset X$}
We will study the existence of terminal contractions by considering cases with respect to the type of singularities 
that the general section $S$ of $X$ through $\Gamma$ has. To do this it will be necessary to obtain normal forms 
for the equation of $\Gamma \subset S \subset X$. 

There is a classification of pairs $(S, \Gamma)$ when $S$ is DuVal and $\Gamma \subset S$ a smooth curve. 
\begin{theorem}[~\cite{Jaf92}]\label{surface_normal_forms}
Let $ 0 \in \Gamma \subset S$ be a smooth curve in a DuVal surface $S$. Then under suitable coordinates, 
$S$ is given by $g=0$ and $\Gamma$ by an ideal $I$ as follows:
\[
\begin{array}{cccc}
\text{Type}             &          g          &          I                                  &     \text{Where}   \\ \hline
A_n           &       z^{n+1}-xy    &    I_k=(x-z^k,y-z^{n+1-k})           &             n \geq 1, \; 1 \leq k \leq \frac{n+1}{2}\\   
D_n           &   x^2 + y^2z-z^{n-1}       &        (x,z)                                 &   n \geq 4 \\
D_{2n}           &        x^2+y^2z-z^{2n-1}  &   (x,y-z^{n-1})         &          n \geq 3 \\
D_{2n+1}         &        x^2+y^2z-z^{2n}    &   (x-z^n,y)             &          n \geq 2 \\    
E_6              &        x^2+y^3-z^4        &   (y,x+z^2)             &                     \\
E_7              &        x^2+y^3+yz^3       &   (x,y)                 &              
\end{array}
\]
\end{theorem}
Let $U \la S$ be the minimal resolution of $S$ and let $E_i$ be the exceptional curves. 
The different forms for the defining ideal of $\Gamma$ correspond to the possible position of $\Gamma$ in the fundamental cycle. 

For the purpose of this paper it will be more convenient to fix the equations of $\Gamma$ and vary this of $S$. The next lemma 
does this and it also shows the relation between the position of $\Gamma$ in the fundamental cycle and its defining ideal. It essentially 
follows from~\cite{Jaf92}. 
\begin{lema}
With assumptions as before. 
\begin{enumerate}
\item Suppose that $0 \in S$ is $A_n$, and that the fundamental cycle is \[
\stackrel{E_1}{\circ} \mbox{\noindent---} \stackrel{E_2}{\circ} \mbox{\noindent ---} 
\cdots  \mbox{\noindent---} \stackrel{E_{n-1}}{\circ} \mbox{\noindent ---} \stackrel{E_n}{\circ}. \]
If $\Gamma$ intersects the $E_k$ then it is given by the ideal $I_k$.
\item Suppose that $0 \in S$ is $D_n$ and that the fundamental cycle is \[
\begin{array}{cc}
                                                                 &  \stackrel{E_{n-1}}{\circ} \\
                                                                 & \mid                       \\
\stackrel{E_1}{\circ} \mbox{\noindent ---} \stackrel{E_2}{\circ} \mbox{\noindent ---}  \cdots \mbox{\noindent ---}  & \stackrel{E_{n-2}}{\circ}  \\
                                                                  & \mid \\
                                                                  & \stackrel{E_n}{\circ}
\end{array}. \]
$\Gamma$ can only intersect $E_1$, $E_{n-1}$ or $E_n$. If it intersects $E_1$ then it is given by $I=(x,z)$ and $S$ by $x^2 + y^2z+z^{n-1}=0$. 
If it intersects $E_{n-1}$ or $E_{n}$, then by an obvious change of variables:
\begin{enumerate}
\item In the $D_{2n}$ case, $S$ is given by $x^2+y^2z+2yz^n=0$ and $\Gamma$ by $I=(x,y)$.
\item In the $D_{2n+1}$ case, $S$ is given by $x^2+y^2z+2xz^n=0$ and $\Gamma$ by $I=(x,y)$.
\end{enumerate} 
\end{enumerate}
\end{lema}
\begin{definition}\label{DF}
Let $\Gamma $ be a smooth curve on a surface $S$. Suppose that $S$ has exactly one singular point on $\Gamma$ 
which is of type $D_n$. Let $f:U \la S$ be its minimal resolution and $\Gamma^{\prime} = f_{\ast}^{-1}\Gamma$. Then:
\begin{enumerate}
\item $\Gamma \subset S$ will be called of type $FD_l$ if $\Gamma^{\prime}$ intersects $E_1$ in the minimal resolution of $S$. 
\item $\Gamma \subset S$ will be called of type $FD_r$ if $\Gamma^{\prime}$ intersects $E_n$ or $E_{n-1}$.
\end{enumerate}
\end{definition}
Next we will derive the simplest possible normal forms for $0 \in \Gamma \subset S \subset X$. To do this it is necessary to obtain some 
properties of $S$. 

First we will show that if $0 \in X$ is $cA_n$ then $0 \in S$ is $A_m$ for the general section $S$ of $X$ containing $\Gamma$.
\begin{lema}
Let $0 \in \Gamma \subset X$, $\Gamma$ a smooth curve and $0 \in X$ a $cA_n$ $3$-fold singularity. 
Then the general hyperplane section $S$ of $X$ containing $\Gamma$ is $A_m$. 
\end{lema}
\begin{proof}
Since $\Gamma$ is smooth, we can assume that it is given by $x=y=z=0$, and $X$ by \[
q(x,y,z)+t[\sum_k \phi_k(x,y,z)t^k]=0\]
and the section given by $t=0$ is $A_n$. Therefore\[
q(x,y,z)=q_2(x,y,z)+\sum_{i\geq 3}q_i(x,y,z)\]
with $q_2(x,y,z)$ a quadratic that is not a square. 

A general hyperplane through $\Gamma$ is $x=by+cz$. Then $S$ is given by \[
q_2(by+cz,y,z)+\sum_{i \geq 3}q_i (by+cz,y,z) + t [\sum_i \phi_i(by+cz,y,z)t^i]=0. \]
It's quadratic term is \[
\psi_2=q_2(by+cz,y,z)+tl(by+cz,y,z)+at^2 \]
where $l(by+cz,y,z)$ is linear. $0 \in S$ is $A_n$ iff $\psi_2$ is not a square. If it is a square, then putting $t=0$ 
it follows that $q_2(by+cz,y,z)$ is also a square. 

\textbf{Claim}: If $q_2(by+cz,y,z)$ is a square for all $b$, $c$ then $q_2(x,y,z)$ is also a square. 

To see this suppose that \[
q_2(x,y,z)=a_1x^2+a_2y^2+a_3z^2+a_4xy+a_5xz+a_6yz.\]
Then
\begin{gather*}
q_2(by+cz,y,z)=\\
=(a_1b^2+a_2+a_4b)y^2+(a_1c^2+a_3+a_5c)z^2+(2bca_1+a_4c+a_5b+a_6)yz. 
\end{gather*}
This is a square iff \[
4(a_1c^2+a_3+a_5c)(a_1b^2+a_2+a_4b)=(2bca_1+a_4c+a_5b+a_6)^2 \;\; \forall b, \; c.\]
Hence $a_1=a_4=a_5=0$ and $4a_2a_3=a_6^2$. But then \[
q_2(x,y,z)=a_2y^2+a_3z^2+a_6yz \]
is a square.
\end{proof}
Now let $f: Y=B_{\Gamma}X \la X$ as before. Let $f^{-1}(\Gamma)=E_1+dE_2$. We want to relate $d$ with some quantity 
on the general hyperplane section $S$ of $X$, through $\Gamma$. 
\begin{lema}
Let $f: Y=B_{\Gamma}X \la X$. Assume that $Y$ is canonical. Suppose that $f^{-1}(\Gamma)=E_1+dE_2$. Then there is a DuVal 
section $S$ of $X$ containing $\Gamma$ such that if $g:S^{\prime}=f_{\ast}^{-1}S: \la S$, then
\[
g^{-1}(\Gamma)=\Gamma^{\prime} + dE 
\]
and $E \cong \pone$ is the $g$-exceptional curve. In fact this is true for the general the general section $S$ of $X$ through $\Gamma$. 
\end{lema}
This way we see that somehow $d$ controls the type of the general hyperplane section through $\Gamma$. This result will be usefull 
later when we try to get normal forms for the equations of $0 \in \Gamma \subset S \subset X$.
\begin{proof}
Let $Q \in L=E_1 \cap E_2$ be a cDV point and let $S^{\prime} \subset Y$ be general through $Q$. Then $Q \in S^{\prime}$ is DuVal. 
Since $S^{\prime}$ is general it has the property that  \[
S^{\prime} \cap (E_1+dE_2)=\Gamma^{\prime} + dE \]
and $\Gamma^{\prime}$ maps to $\Gamma$. Moreover, for a general $\delta \subset E_1$, \[
\delta \cdot S^{\prime} =1 \]
and hence \[
\delta \cdot ( S^{\prime} + E_1+dE_2)=0.\]
Therefore, there is a Cartier divisor $S$, $ 0 \in \Gamma \subset S \subset X$ such that \[
f^{\ast}S = S^{\prime} + E_1 + dE_2. \]
By adjunction \[
K_{S^{\prime}}=(K_Y+S^{\prime})\mid_{S^{\prime}}=f^{\ast}(K_X +S) \mid_{S^{\prime}}=f^{\ast}K_S. \]
Therefore $S$ is canonical with the required property.
\end{proof}
The next proposition relates $d$ with the type of singularities of $S$. 
\begin{proposition}
Let $P \in \Gamma \subset S$. $\Gamma$ a smooth curve and $P \in S$ DuVal. Let $g: S^{\prime} \la S$ be the blow up 
of $S$ along $\Gamma$. Let $E$ be the $g$-exceptional curve which is necessarily irreducible. Suppose that 
$g^{-1}(\Gamma)=\Gamma^{\prime}+dE$. Let $f:U \la S$ be the minimal resolution of $S$, $E_i$ the $f$-exceptional curves 
and $\Gamma^{\prime \prime}$ the birational transform of $\Gamma$ in $U$. Then
\begin{enumerate}
\item Suppose that $P \in S$ is $A_n$ and that $\Gamma^{\prime\prime}$ intersects $E_k$, $1 \leq k \leq [\frac{\text{n+1}}{2}]$. Then 
\begin{enumerate}
\item \[
[f^{-1}(\Gamma)]=\Gamma^{\prime\prime} +\sum_{i=1}^{k-1}iE_i + k \sum_{i=k}^{n-k+1}E_i + \sum_{i=1}^{k-1}(k-i)E_{n-k+i+1}.
\]
\item $E=E_{n-k+1}$, and $d=k$.
\end{enumerate}
\item Suppose that $P \in S$ is $D_n$. Then $\Gamma^{\prime\prime}$ intersects one of the edges of the dual graph of $S$. Then
\begin{enumerate}
\item With notation as in lemma 4.2.2, if $\Gamma^{\prime\prime}$ intersects $E_1$, then $E=E_1$ and $d=2$. Moreover $S^{\prime}$ 
is smooth along $\Gamma^{\prime}$ and 
\[
[f^{-1}(\Gamma)]=\Gamma^{\prime\prime} +2\sum_{i=1}^{n-2}E_{i} +E_{n-1}+E_{n} \]
\item If $\Gamma^{\prime\prime}$ intersects $E_{n-1}$ then 
%\[
%E=\left\{
%\begin{array}{cc}
%E_{n-1}     & \text{if n is even, and d=n/2. $S^{\prime}$ is smooth along $\Gamma^{\prime}$.} \\
%E_n         & \text{if n is odd, and $d=\frac{n-1}{2}$}
%\end{array} \right.
%\]
\begin{enumerate}
\item If $n$ is even, then $E=E_{n-1}$, $d=n/2$, $S^{\prime}$ is smooth along $\Gamma^{\prime}$ and \[
[f^{-1}(\Gamma)]=\Gamma^{\prime\prime}+E_1+2\sum_{i=2}^{n-2}E_i +\frac{n}{2}E_{n-1}+\frac{n-2}{2}E_n \]
\item If $n$ is odd, then $E=E_n$, $d=\frac{n-1}{2}$, and \[
[f^{-1}(\Gamma)]=\Gamma^{\prime\prime}+E_1+2\sum_{i=2}^{n-2}E_i +\frac{n-1}{2}E_{n-1}+\frac{n-1}{2}E_n.
\]
\end{enumerate}
\end{enumerate} 
\end{enumerate}
\end{proposition}
\begin{corollary}
Let $P \in \Gamma \subset S \subset X$, $\Gamma$ smooth and $P \in X$ a cDV point. 
Let $f:Y=B_{\Gamma}X \la X$, $f^{-1}(\Gamma)=E_1+dE_2$ and $L=E_1 \cap E_2$. 
If $P \in S$ is a $D_{2n}$ singularity, then $Y$ has finitely many singular points on $L$ and therefore 
a $\mathbb{Q}$-factorialization of $E_1$ is obtained by blowing up $E_1$.
\end{corollary}
\begin{proof}
From the previous proposition it follows that $S^{\prime}=f_{\ast}^{-1}S$ is smooth at one point of $L$. Therefore since $S^{\prime}$ 
is Cartier, $Y$ is also smooth at this point and hence has finitely many singularities along $L$.
\end{proof}
\begin{proof}[Proof proposition 5.6]
Let $Z=[f^{-1}(\Gamma)]$. This is an integral cycle. From the properties of the blow up it follows that 
$g^{-1}(\Gamma)$ is Cartier and hence $Z \cdot E =-1$ and $Z \cdot E_i =0$ for all $i$ such that $E_i \neq E$. 
Moreover, $d$ is just the coefficient of $E$ in $Z$. 

We will only do the case  when $0 \in S$ is $D_n$, with $n$ odd, and $\Gamma^{\prime\prime}$ intersects $E_{n-1}$. 
The rest is similar. In fact the $A_n$ is simpler. 

Observe that since $S$ has embedding dimension 4 and $S^{\prime}$ is the blow up of a smooth curve, it follows that $g^{-1}(0)=E=\pone$. 
Therefore $E$ appears with coefficient 1 in the fundamental cycle and hence it must be one of the edges. 

So let $Z=\Gamma^{\prime\prime}+\sum_{i}a_i E_i $, $a_i \in \mathbb{N}$.

There are three cases to be considered. Only one will give an integral cycle and this will be the answer. 

\textbf{Case 1}: Check if $E=E_{n}$. The relations $Z \cdot E_{n} =-1$ and $Z \cdot E_i =0$, $\forall \; i \neq n$ give the system 
of equations 
\begin{gather*}
-2a_1+a_2=0 \\
a_1-2a_2+a_3=0 \\
\vdots \\
a_{n-3}-2a_{n-2}+a_{n-1}+a_n=0 \\
a_{n-2}-2a_{n-1}+1=0\\
a_{n-2}-2a_n =-1
\end{gather*}
It is easy to see that the solution of this system is $a_i=i$ for $2\leq i \leq n-2$, and $a_{n-1}=a_{n}=(n-1)/2$. This solution 
gives an integral cycle. Similarly we see that the cases $E=E_{n-1}$ and $E=E_1$ do not give integer cycles and hence are not possible.

Therefore, \[
[f^{-1}(\Gamma)]=\Gamma^{\prime\prime}+E_1+2\sum_{i=2}^{n-2}E_i +\frac{n-1}{2}E_{n-1}+\frac{n-1}{2}E_n, \]
$E=E_n$ and $d=(n-1)/2$.
\end{proof}
We are now in position to get normal forms for the equations $0 \in \Gamma \subset S \subset X$ in the case 
that $0 \in S$ is $D_n$ and $S$ is general through $\Gamma$. 
We will not treat the general $A_n$ case here and so I will not attempt to write normal forms in this case. 
However, normal forms for the case when $0 \in S$ is $A_3$ will be given in the proof of Theorem~\ref{An}. 

\begin{proposition}\label{normalforms}
Suppose $0 \in \Gamma \subset S \subset X$, $S$ general. Suppose that $0 \in S$ is $D_n$. 
Then under suitable choice of coordinates, $0 \in \Gamma \subset S \subset X$ is given by:
\begin{enumerate}
\item If $\Gamma \subset S$ is of type $FD_l$, then 
\[
\begin{array}{cc}
X: \; x^2+y^2z+z^{n-1}+t\phi_{\geq 2}(y,z,t)=0,  &  I_{\Gamma} = (x,z,t). 
\end{array} \]
Moreover no $y^k$ appears in $\phi_{\geq 2}(y,z,t)$ for any $k$. For $n=4$ this is the only possibility. \\
If $n \geq 5$ then $\phi_{2}=0$. That is $X$ is given by \[
x^2+y^2z+z^{n-1}+t\phi_{\geq 3}(y,z,t)=0 \] 
and again there is no $y^k$ appears in $\phi_{\geq 3}(y,z,t)$. 
\item Suppose that $\Gamma \subset S$ is of type $FD_r$. Then 
\begin{enumerate}
\item If $n$ is even then $0 \in \Gamma \subset S \subset X$ is given by:
\[ \begin{array}{cc}
x^2+y^2z+2yz^{n/2} + t \phi_{\geq 2}(y,z,t)=0, &  I_{\Gamma} = (x,y,t). 
\end{array} \]
Moreover, There is no $z^k$ in $\phi_{\geq 2}(y,z,t)$ for any $k$. 
\item If $n$ is odd then $0 \in \Gamma \subset S \subset X$ is given by:
\begin{enumerate}
\item \[ \begin{array}{ccc}
x^2 +y^2z+2xz^{\frac{n-1}{2}} + t [ axz^k + \phi_{\geq 2}(y,z,t)]=0, &k \geq 1 & I_{\Gamma} = (x,y,t).
\end{array} \]
No $yz$ or $z^{\nu}$ appear in $\phi_{\geq 2}(y,z,t)$ for any $\nu$. $a \in \mathbb{Z}$. \\
\item Alternatively, the equation can take the form \[
x^2 +y^2z+2xz^{\frac{n-1}{2}} + t [xz \psi (z,t) + axt^{\frac{n-3}{2}}+bxt^k+\phi_{\geq 2}(y,z,t)]=0, 
\]
$k \geq 1$, $ I_{\Gamma} = (x,y,t)$ and in this case, $y^2$, $yz$, or $z^{\nu}$ do not appear in $\phi_{\geq 2}(y,z,t)$, for any $\nu$.
\end{enumerate}
\end{enumerate}
\end{enumerate}
\end{proposition}
Sometimes it is better to have $2.b$ and sometimes $2.a$. The existence of $y^2$ may complicate calculations. 
\begin{proof}
We will apply the following methods:
\begin{enumerate}
\item The Weierstrass preparation theorem
\item The elimination of the $y^{n-1}$ -term from the polynomial $a_ny^n+a_{n-1}y^{n-1}+\cdots $ by a coordinate change 
$y \la y-a_{n-1}/na_n $ when $a_n$ is a unit.
\item Let $M_1$, $M_2$, $M_3$, $M_4$ be multiplicatively independent monomials in the variables $x$, $y$, $z$, $t$. 
Then any power series of the form $M_1 \cdot (unit) + M_2 \cdot (unit) + M_3 \cdot (unit) + M_4 \cdot (unit)$ is 
equivalent to $M_1+M_2+M_3+M_4$ by a suitable coordinate change $x \mapsto x \cdot (unit)$, $y \mapsto y \cdot (unit)$, 
$z \mapsto z \cdot (unit)$, $t \mapsto t \cdot (unit)$. 
\end{enumerate}

\textbf{Case 1}: Suppose that $\Gamma^{\prime\prime}$ intersects $E_1$ in the dual graph. 

Then by lemma 4.2, $0 \in \Gamma \subset X$ is given by \[
x^2+y^2z+z^{n-1}+t\phi_{\geq 1}(x,y,z,t)=0.\]
Apply the Weierstrass preparation theorem on $x^2$ to eliminate $x$ from $\phi$. Moreover, by lemma 4.4, if $\phi$ has linear terms, 
then $0 \in X$ is $cA_n$ and therefore the general section $\Gamma \subset S \subset X$ is $A_m$ which is not possible by our 
assumptions. Therefore $\Gamma \subset S \subset X$ is given by \[
x^2+y^2z+z^{n-1}+t\phi_{\geq 2}(y,z,t)=0, \]
and $\Gamma =(x,z,t)$.

Now suppose that \[
\phi_{\geq 2}(y,z,t)=f_{\geq 2}(y) + \Phi_{\geq 2}(y,z,t), \]
with no $y^k$ in $\Phi_{\geq 2}(y,z,t)$. Then write $f(y)=y^k \cdot (unit)$, $k \geq 2$. Then the equation of $X$ is \[
x^2+y^2 (z+ty^{k-2} \cdot (unit))+z^{n-1} +t \Phi_{\geq 2}(y,z,t)=0. \]
The change of variables $z \mapsto z-ty^{k-2} \cdot (unit)$ will give the normal form claimed by the first part of $1.$

We will get further restrictions on $\phi$ in the case $n \geq 5$ later. 

\textbf{Case 2}: Suppose that $\Gamma^{\prime\prime}$ intersects $E_{n-1}$ or $E_n$ and that $n=2m$. 

Again by lemma 4.2, and by applying the Weierstrass preparation theorem it follows 
that $\Gamma \subset S \subset X$ is given by \[
x^2 +y^2z+2yz^m+t\phi_{\geq 2}(y,z,t)=0, \]
and $I_{\Gamma} =(x,y,t)$. We want to eliminate $z^k$ as well. Suppose as before that \[
\phi_{\geq 2}(y,z,t)=f_{\geq 2}(z) + \Phi_{\geq 2}(y,z,t)= z^k \cdot (unit) +  \Phi_{\geq 2}(y,z,t). \]
\textbf{Claim}: $k \geq m$. 

If this is true, then the equation of $X$ can be written as \[
x^2+y^2z+2z^m (y+tz^{k-m}\cdot (unit))+t \psi_{\geq 2}(y,z,t)=0.\]
Then the change of coordinates $y \mapsto y-tz^{k-m}\cdot (unit)$ gives the claimed normal form. 

To see the claim, we will show by using lemma 4.5 that if $k <m$ then there is a $D_{s}$ section $\Gamma \subset T \subset X$ with 
$s < n$. By that lemma, we just have to show that if $f: Y=B_{\Gamma}X \la X$ and $f^{-1}(\Gamma)=E_1+dE_2$, then $d<m$. 

Calculate the blow up. In the chart $x=xt$, $y=yt$, $Y$ is given by \[
x^2t+y^2tz+2yz^m+z^k \cdot (unit) + \Phi_{\geq 2}(yt,z,t). \]
For $t=0$ we get that $I_{E_1} = (t, 2yz^{m-k}+ (unit))$, $I_{E_2}=(t,z)$ and \[
f^{-1}(\Gamma)=E_1+kE_2 
\]
and $k < m$. So we obtain $2.a$.

\textbf{Case 3}: As in case 2 but $n=2m+1$. 

Again by the preparation theorem the equation of $0 \in \Gamma \subset S \subset X$ takes the form 
\begin{equation}
x^2 +y^2z+2xz^{m} + t\phi_{\geq 2}(y,z,t)=0 
\end{equation}
with $I_{\Gamma} = (x,y,t)$. 

We want to eliminate $z^k$. As before, write \[
\phi_{\geq 2}(y,z,t)=f_{\geq 2}(z) + \Phi_{\geq 2}(y,z,t)= z^k \cdot (unit) +  \Phi_{\geq 2}(y,z,t), \]
and for the same reasons as in the previous case it follows that $k \geq m$. Moreover, in this case, $k > m$. 
Because suppose that $k=m$. Then by using $(1)$, one can calculate, just like the proof of Theorem~\ref{Dn} that if 
$f:Y \la X$ is the blow up of $X$ along $\Gamma$, and $f^{-1}(\Gamma)=E_1+mE_2$, then $Y$ has finitely many singularities 
along $L=E_1 \cap E_2$. But this is not possible as follows from lemma~\ref{diff}.

Then $(1)$ becomes \[
x^2+y^2z+2z^m(x+tz^{k-m}\cdot (unit)) +t \Phi_{\geq 2}(y,z,t)=0. \]
The change of variables $x \mapsto x-tz^{k-m}\cdot (unit)$ and $t \mapsto t/(unit)$ makes the equation of $X$ \[
x^2+y^2z+2xz^m+t[axz^{\nu}+\phi_{\geq 2}(y,z,t)]=0,\]
$\nu \geq 1$ and $I_{\Gamma} = (x,y,t)$ and no power $z^k$ appears in $\phi_{\geq 2}(y,z,t)$. 

To eliminate $yz$, write $\phi_{\geq 2}(y,z,t)=byz + \{other\}$ and make the change of coordinates $y \mapsto y-bt/2$. This will 
produce the required normal form. 

Sometimes it is preferable to eliminate $y^2$ too. To do that write $\phi_{\geq 2}(y,z,t)=cy^2+\{other\}$ and make the change 
of variables $z \mapsto z-bt$. This will produce the normal form in $2.b.ii.$

The only thing left to show is the second statement of $1.$ for $n \geq 5$. The point is that 

\textbf{Claim}: If $\phi_2(y,z,t)\neq 0$ then there is a $D_4$ section $S$ of $X$ through $\Gamma$.

We already know that in the $DF_l$ case it is possible to write the equations of $\Gamma \subset S \subset X$ as\[
G=x^2+y^2z+z^{n-1} +t \phi_{\geq 2}(y,z,t)=0, \]
and no power of $y$, $y^k$, appears in $\phi_{\geq 2}(y,z,t)$. Let \[
\phi_{\geq 2}(y,z,t)=\phi_2(y,z,t)+\phi_{\geq 3}(y,z,t). \]
Let $\phi_2(y,z,t)=a_1z^2+a_2t^2+a_3yz+a_4yt+a_5zt$. Let $S$ be the section given by $t=\lambda z$. This is given by \[
F(x,y,z)=x^2+y^2z+z^{n-1}+\lambda z \phi_2(y,z,\lambda z)+\lambda z \phi_{\geq 3}(y,z,\lambda z)=0, \]
and \[
\phi_2(y,z, \lambda z)=(a_1+a_2 \lambda^2+a_5 \lambda)z^2+(a_5+a_4 \lambda )yz.\]
If $\phi_2 \neq 0$, then for general $\lambda$ at least one of the coeffiecients of $z^2$ or $yz$ in not zero. 

\textbf{Case 1:} Both coefficients are nonzero. i.e., \[
\begin{array}{ccc}
a=a_1+a_2\lambda^2+a_5 \lambda \neq 0 & \text{and} & b=a_3+a_4 \lambda \neq 0.
\end{array} \]
Hence $\phi_2(y,z,\lambda z)=az^2+byz$, $ab \neq 0$ and \[
F=x^2+y^2z+z^{n-1}+\lambda z (az^2+byz)+ \lambda z \phi_{\geq 3}(y,z,\lambda z).\]
Now write \[
\phi_{\geq 3}(y,z,\lambda z)=zf_{\geq 2}(y)+z^2f_{\geq 1}(y,z,\lambda z).\]
Also write $f_{\geq 2}(y)=y^k f(y)$, $k \geq 2$. Then 
\begin{gather*}
F=x^2+y^2z+yz^2[\lambda b + \lambda y^{k-1}f(y)]+z^3[\lambda a +z^{n-4}+ \lambda f_{\geq 1}(y,z,\lambda z)]=\\
=x^2+y^2z+yz^2 \cdot u_1 +z^3 \cdot u_2
\end{gather*}
where, \[
\begin{array}{ccc}
u_1=\lambda b + \lambda y^{k-1}f(y)   &  \text{and} & u_2=\lambda a +z^{n-4}+ \lambda f_{\geq 1}(y,z,\lambda z)
\end{array} \]
are units. Hence \[
F=x^2+z[y^2+yz \cdot u_1 +z^2 \cdot u_2]. \]
Now eliminate the $yz$-term by the change of variables $y \mapsto y-zu_1/2$. Hence the equation becomes \[
F=x^2+z[y^2-\frac{\textstyle z^2u_1^2}{\textstyle 4}+z^2u_2]=x^2+y^2z+z^3(u_2-\frac{\textstyle u_1^2}{\textstyle 4}). \]

We will now consider cases with respect to the nature of $\delta = u_2-\frac{\textstyle u_1^2}{\textstyle 4}$. 

Suppose that $\delta$ is a unit. Then \[
 F \sim x^2+y^2 z+z^3 \]
and therefore $S$ is $D_4$. 

Suppose now that $\delta$ is not a unit. First check when this happens. By looking at how $u_1$, $u_2$ are defined, we see 
that \[
\delta = (\lambda a - \frac{\textstyle \lambda^2 b^2}{\textstyle 4}) + \{\text{higher}\}.\]
Therefore $\delta$ is not a unit iff\[
a-\frac{\textstyle \lambda b^2}{\textstyle 4}=0, \;\; \forall \lambda.\]
But this implies that \[
a_1+\lambda^2 a_2+\lambda a_5=\frac{\textstyle 1}{\textstyle 4} \lambda (a_3+a_4\lambda)^2 \]
for all $\lambda$, and hence \[
a_1=a_2=a_4=4a_5-a_3^2=0.\]
Therefore \[
\phi_2(y,z,t)=\frac{\textstyle 1}{\textstyle 4} a^2 zt +ayz.\]
Then 
\begin{gather*}
G=x^2+y^2z+z^{n-1}+t(\frac{\textstyle 1}{\textstyle 4} a^2 zt+ayz)+t\phi_{\geq 3}(y,z,t)=\\
=x^2+z(y^2+\frac{\textstyle 1}{\textstyle 4}a^2t^2+ayt)+z^{n-1}+t\phi_{\geq 3}(y,z,t)=\\
=x^2+z(y+\frac{\textstyle 1}{\textstyle 2}at)^2+z^{n-1}+t\phi_{\geq 3}(y,z,t).
\end{gather*}
Make the change of variables $y \mapsto y-\frac{1}{2}at$ to bring the equation in the form \[
x^2+y^2z+z^{n-1}+t\phi_{\geq 3}(y,z,t)=0, \]
which is what we want. 

We must now check what happens if $a=0$ or $b=0$ for all $\lambda$. 

Suppose that $a=0$. Hence $a_1=a_2=a_5=0$ and therefore \[
\phi_2(y,z,t)=a_3yz+a_4yt.\]
Then \[
\phi_2(y,z,\lambda z)=(a_3+a_4 \lambda)yz \]
which is nonzero for general $\lambda$. Let $S$ again be the section defined by $t=\lambda z$. This is given by 
\begin{gather*}
F=x^2+y^2z+z^{n-1}+\lambda z \phi_2 (y,z,\lambda z)+\lambda z \phi_{\geq 3}(y,z,\lambda z)=\\
x^2+y^2z+z^{n-1}+\lambda z^2y(a_3+a_4\lambda)+\lambda z \phi_{\geq 3}(y,z,\lambda z)=\\
x^2+z(y^2+\lambda ( a_3+a_4 \lambda)yz)+z^{n-1}+\lambda z\phi_{\geq 3}(y,z,\lambda z)=0.
\end{gather*}
Now eliminate the $yz$ term by the change of variables \[
y \mapsto y-\lambda (a_3+a_4 \lambda)z/2\]
to get \[
F=x^2+z(y^2-\gamma^2z^2)+z^{n-1}+\lambda z \phi_{\geq 3}(y,z,\lambda z)=0 \]
where $\gamma = \lambda (a_3+a_4 \lambda)z/2$. Since there is no isolated power $y^k$ it is possible to write \[
\phi_{\geq 3}(y,z,\lambda z)=zf_{\geq 2}(y)+z^2f_{\geq 1}(y,z)\]
and $f_{\geq 2}(y)=y^kf(y)$ for $k \geq 2$. Hence 
\begin{gather*}
F=x^2+y^2z-\gamma^2z^3+z^{n-1}+\lambda z^2y^kf(y)+\lambda z^3 f_{\geq 1}(y,z)=\\
=x^2 + y^2z[1+\lambda z y^{k-2} f(y)]+z^3[-\gamma + \text{\{higher\}}]=\\
=x^2 +y^2z \cdot (unit) +z^3 \cdot (unit) \sim x^2+y^2z+z^3,
\end{gather*}
and hence $S$ is $D_4$ which is not possible.

The last case to consider is when $a_3=a_4=0$. In this case, \[
\phi_2(y,z,t)=a_1z^2+a_2t^2+a_5zt,\]
and therefore \[
\phi_2(y,z,\lambda z)=\gamma z^2 \]
where $\gamma=a_1+a_2 \lambda^2 + a_5 \lambda$. Hence \[
F=x^2+y^2z+z^{n-1}+\lambda z (\gamma z^2 + \phi_{\geq 3}(y,z,\lambda z)). \]
Again write $\phi_{\geq 3}=zy^kf(y)+z^2f_{\geq 1}(y,z)$, and hence 
\begin{gather*}
F=x^2+y^2z+z^{n-1}+\lambda z[\gamma z^2 + zy^kf(y)+z^2f_{\geq 1}(y,z)]=\\
=x^2+y^2z+z^{n-1}+\gamma \lambda z^3 + \lambda z^2 y^kf(y)+\lambda z^3f_{\geq 1}(y,z)=\\
=x^2+y^2z[1+\lambda z f(y)]+z^3[\gamma \lambda  + z^{n-3}+\lambda f_{\geq 1}(y,z)]=\\
x^2+y^2z \cdot (unit) +z^3 \cdot (unit) \sim x^2+y^2z+z^3,
\end{gather*}
and therefore $S$ is again $D_4$ which is impossible.

This concludes the proof of the proposition.
\end{proof}

%% file: section3.tex
\section{The $A_1$, $A_2$, and $A_3$ cases.}
In this section we will study the existence of terminal contractions in the case when $0 \in S$ is $A_1$, $A_2$, $A_3$ and 
$\Gamma^{\prime\prime}$ intersects the edge of the dual graph. 
\begin{theorem}\label{An}
Let $P \in \Gamma \subset S \subset X$. Suppose that $0 \in S$ is $A_n$. Then 
\begin{enumerate}
\item Suppose that $\Gamma^{\prime\prime}$ intersects one edge of the dual graph. 
(In particular, this is always the case when $n=1,\; 2$). Then there is always 
a terminal contraction $E \subset W \la \Gamma \subset X$ contracting a surface $E$ to $\Gamma$. Let $i(W)$ be the index of $W$. Then \[
n+1 \leq i(W) \leq 2n.\] 
In particular if $n=1$ then $i(W)=2$ and the index 2 points of $W$ are $cA$. If $n=2$ then $i(W) \in \{3,4\}$. 
Moreover, \[
R(\Gamma, X)=\oplus_{d \geq 0}\I_{\Gamma, X}^{(d)}\]
is finitely generated by elements of degrees at most $2n$.
\item Suppose that $n=3$. Then
\begin{enumerate}
\item If $\Gamma^{\prime\prime}$ intersects one edge of the dual graph, then there is a terminal 
contraction as follows from the previous part.
\item If $\Gamma^{\prime\prime}$ intersects $E_2$, i.e. the middle of the dual graph, then write the equation of 
$\Gamma \subset S \subset X$ as
\[x^2+y^2+f_{\leq 3}(y,z,t)+f_{\geq 4}(y,z,t)=0,\]
no $y^2$ appears in $f_{\leq 3}(y,z,t)$ and $I_{\Gamma} =(x,y,t)$.( This is always possible). 
Then a terminal contraction exists iff there is no proper 
irreducible component of $(f_{\leq 3}(y,z,t)=0) \subset \mathbb{C}^3$ that goes through the origin. 
If it exists, then it has index 2 and its index 2 points are $cA$. 
\end{enumerate}
\end{enumerate}  
\end{theorem}
\begin{proof}
Fix notation as in the proof of Theorem~\ref{canonicalcase}.

\textbf{Case 1:} Assume that $\Gamma^{\prime\prime}$ intersects one edge of the dual graph. 

Now proceed as in the proof of Theorem~\ref{canonicalcase}. To determine whether or not the contraction is terminal, 
it suffices to check whether $Z$ has isolated terminal singularities away from $E_2^Z$, or not. 
From proposition 4.6 it follows that $d=1$ and therefore \[
f^{-1}(\Gamma)=E_1+E_2. \]
$E_1$, $E_2$ are both smooth and hence $Y=B_{\Gamma}X$ is $cA_{\ast}$ along $L=E_1 \cap E_2$. Let $C \subset Z$ be a $g$-exceptional 
curve. It must lie over a cDV point and therefore $Z$ can have at most finitely many terminal singularities along $C$. 
Hence $W$ is terminal. 

Now to find the index of $W$. Let $b \in \mathbb{N}$ such that \[
K_Z =p^{\ast}K_W+bE_2^Z. \]
\textbf{Claim:} $E_2^Z$ has index $n$. 

Let $S^{\prime}=f^{\ast}S$. Then $S^{\prime}$ has exactly 1 singular point which is $A_{n-1}$. This follows since $\Gamma^{\prime\prime}$ 
intersects the edge of the dual graph. So at the generic point of $L$, $Y$ has a $A_{n-1}$ point as follows from 
lemma 4.5 and proposition 4.6. At this point, $E_1$, $E_2$ correspond to two lines at the edge of the dual graph
Therefore $nE_1$, $nE_2$ are Cartier at all but finitely many points. 

Since $E_1+E_2$ is Cartier, it follows that the singularities of $Y$ lie on $L=E_1 \cap E_2$. In fact if $n \geq 2$, then $Y$ 
is singular along $L$. If not then by lemma 4.5, there is a $\Gamma \subset S \subset X$ such that $S^{\prime}=B_{\Gamma}S$ is smooth 
along $S^{\prime} \cap E_1$. But by proposition 4.6 this implies that $S$ must be $A_1$ which is not possible. 

Moreover, since $Z$ has index 1, it has hypersurface singularities and therefore by lemma 2.2 the indices 
of $E_1^Z$ and $E_2^Z$ can be computed at any point of $L^Z$. Therefore, they have index $n$. 

Let $l \subset E_2^{Z}$ be the birational trasform of a line in $E_2 = \mathbb{P}^2$ contracted by $p$. 
Then, $l \cdot E_2^Z = - a/n$, $a \in \mathbb{N}$. Moreover, \[
l \cdot E_2^Z + l \cdot E_1^Z =-1 \] 
and clearly $l \cdot E_1^Z \leq 1$. Hence \[
-2 \leq l \cdot E_2^Z \leq -1.\]
Therefore $n \leq a \leq 2n$. Moreover, $K_Z \cdot l =K_Y \cdot l =-1$ and hence $b=n/a$. In fact $a=n$ is not possible. If 
it was then \[
K_Z = p^{\ast}K_W + E_2^Z.\]
Moreover, \[
\begin{array}{ccc}
K_W=q^{\ast}K_X+E_1^W  & \text{and}  &     K_Z=h^{\ast}K_X+E_1^Z+E_2^Z 
\end{array}\]
where $h=f \circ g: Z\la X$. Combining the above relations it follows that \[
p^{\ast}E_1^W \cong E_1^Z \]
which is not possible. Now the claim about the indices follows immediately. 

The statement about the type of singularities follows from~\cite[Theorem 4.7]{Ko-Mo92}.

\textbf{Case 2:} $n=3$. The only case to study is when $\Gamma^{\prime\prime}$ intersects the middle of the dual graph. 
It will be necessary to obtain a normal form for the equation $\Gamma \subset S \subset X$. 
\begin{lema}
Let $0 \in \Gamma \subset S \subset X$. Suppose that $0 \in S$ is an $A_3$ singular point and that $\Gamma^{\prime\prime}$ 
intersects the middle of the dual graph in the minimal resolution of $S$. Then under suitable coordinates, 
$0 \in \Gamma \subset S \subset X$ can be written as \[
x^2+y^2+2xz^2+t \phi_{\geq 1}(x,z,t)=0 \]
and $I_{\Gamma}=(x,y,t)$ and no power $z^k$ appears in $\phi_{\geq 1}(x,z,t)$.
\end{lema}
\begin{proof}
By Theorem 4.1, under suitable coordinates it is possible to write $S: xy-z^5=0$ and $\Gamma : x-z^2=y-z^2=0$. The change of 
coordinates $x \mapsto x+z^2$, $y \mapsto y+z^2$ brings it to \[
xy+xz^2+yz^2=0\]
and $\Gamma : x=y=0$. Now let $x \mapsto x-y$, $y \mapsto x+y$ and apply the Weierstrass preparation theorem to $y^2$ to get \[
S: x^2+y^2+2xz^2+t\phi_{\geq 1}(x,z,t)=0\]
and $\Gamma: x=y=t=0$. To eliminate $z^k$ we must show that $z$ does not appear in $\phi_{\geq 1}(x,z,t)$. 
If it does then it is easy to see that 
the general hyperplane section $y=at+bx$ through $\Gamma$ will be $A_1$ which is impossible.

So now proceed as in proposition~\ref{normalforms} to eliminate powers of $z$. 
\end{proof}
The statement about existence of a terminal contraction is proved in exactly the same way as Theorem~\ref{Dn} and I omit its proof.
\end{proof}
\begin{example}
Let $X$ be given by \[
x^2+y^2+2xz^2+t^m=0 \]
and $\Gamma: x=y=t=0$. Then there is no 3-fold terminal contraction, contracting a surface to $\Gamma$.
\end{example}
\begin{proof}
In this case, the section $S: (t=0)$ is an $A_3$ type and the curve intersects the middle of the minimal resolution of $S$. 
Moreover, $f_{\leq 3}(x,z,t)=2xz^2$ and there is a plane through the origin in $f_{\leq 3}(x,z,t)=0$. Therefore by the previous theorem 
there is no terminal contraction.
\end{proof}

%% file: section4.tex
\section{The $D_{2n}$ cases. }
In this section we will study the existence of terminal contractions when $0 \in S$ is a $D_{2n}$ type for general 
$0 \in \Gamma \subset S \subset X$. 
\begin{theorem}\label{Dn}
Let $P \in \Gamma \subset S \subset X$. Suppose that $P \in S$ is a $D_n$ type singular point. Then
\begin{enumerate}
\item Suppose that $n=4$. Write the equation of $P \in \Gamma \subset S \subset X$ as \[
x^2 +f_3(y,z,t)+f_{\geq 4}(y,z,t)=0\]
where $I_{\Gamma} =(x,z,t)$ and $f_3(y,z,t)$ is homogeneous of degree $3$. This is always possible. Then
\begin{enumerate}
\item If $f_3(y,z,t)$  is an irreducible homogeneous cubic, then there is a terminal contraction $W \la X$ of a surface to $\Gamma$. 
$W$ has index 2 and has exactly one singularity which is of $cA_x$ type. Moreover, \[
R(\Gamma, X)=\oplus_{d \geq 0}\I_{\Gamma, X}^{(d)}\]
is finitely generated by elements of degrees 1 and 2.
\item If $f_3(y,z,t)$ is reducible or $0$, then there is no terminal contraction.
\end{enumerate}
\item Suppose $n \geq 5$.
\begin{enumerate}
\item Suppose that $\Gamma \subset S$ is of type $FD_l$. Then there is no terminal contraction. 
\item Suppose that $\Gamma \subset S$ is of type $FD_r$ and $n$ is even. Then write the equation of 
$P \in \Gamma \subset S \subset X$
\[
x^2 +f_3(y,z,t)+f_{\geq 4}(y,z,t)=0\]
with $I_{\Gamma} =(x,y,t)$.
\begin{enumerate}
\item If $f_3(y,z,t)$  is an irreducible homogeneous cubic, then there is a terminal contraction $W \la X$ of a surface to $\Gamma$. 
$W$ has index 2 and has exactly one singularity which is of $cD$ type. Moreover, \[
R(\Gamma, X)=\oplus_{d \geq 0}\I_{\Gamma, X}^{(d)}\]
is finitely generated by elements of degrees 1 and 2.
\item If $f_3(y,z,t)$ is reducible or $0$, then there is no terminal contraction. 
\end{enumerate}
\end{enumerate}
\end{enumerate}
\end{theorem}
\begin{proof}
we will only do the second part of the theorem. The first one is proved in exactly the same way. 

\textbf{Case 1:} $\Gamma^{\prime\prime}$ intersects $E_1$ in the dual graph and $n \geq 5$. 

By proposition~\ref{normalforms}, under suitable coordinates $\Gamma \subset S \subset X$ is given by \[
x^2+y^2z+z^{n-1}+t\phi_{\geq 3}(y,z,t)=0 \]
$\Gamma: x=z=t=0$, and no power $y^k$ appears. Let \[
f : Y=B_{\Gamma}X \la Y\]
be the blow up of $X$ along $\Gamma$. In the chart $x=xt$, $z=zt$ it is given by \[
x^2t+y^2z+z^{n-1}t^{n-2}+\phi_{\geq 3}(y,zt,t)=0.\]
For $t=0$ we see that \[
f^{-1}(\Gamma)=E_1+2E_2 \]
and $E_1: z=t=0$, $E_2: y=t=0$. Moreover $L=E_1 \cap E_2: y=z=t=0$ and it is easy to see that $Y$ has only 
one singular point on $L$, the origin. However it is singular along the line $l: x=y=z=0$ which lies in $E_2$. 

Therefore \[
g : Z=B_{E_1}Y \la Y \]
is a $\mathbb{Q}$-factorialization of $E_1$. In the chart $z=zt$ it is given by \[
x^2+y^2z+z^{n-1}t^{2n-4}+ \frac{1}{t} \phi_{\geq 3}(y,zt^2,t)=0.\]
Write \[
\phi_{\geq 3}(y,zt^2,t)=tf_{\geq 2}(y)+t^2f_{\geq 1}(y,z,t).\]
It is easy to see that $f_{\geq 1}(0,z,0)=0$. Then $Y$ is given by \[
x^2+y^2z+z^{n-1}t^{2n-4}+f_{\geq 2}(y)+tf_{\geq 1}(y,z,t)=0.\]
Now let $C=g^{-1}(0): x=y=t=0$. Then, since $f_{\geq 1}(0,z,0)=0$ it is clear that $Z$ is singular along $C$. Moreover since 
$E_2^{Z} \cong E_2 \cong \p^2$, it is clear that $C$ is not contained in $E_2^{Z}$. 

Now contract $E_2^{Z} \cong \p^2$ over $X$ to get a contraction \[
p: W \la X.\]
However $W$ is only canonical and not terminal since it is singular along $C$. Therefore in this case there is no terminal 
contraction.

\textbf{Case 2:} $\Gamma^{\prime\prime}$ intersects the other edge of the dual graph and $n=2m$ is even. 

Then again by proposition~\ref{normalforms}, under suitable coordinates $\Gamma \subset S \subset X$ is given by \[
x^2+y^2z+2yz^m+t\phi_{\geq 2}(y,z,t)=0,\]
$\Gamma:x=y=t=0$ and no power $z^k$ appears. 
The method is the same as before. Blow up $\Gamma$. In the chart $x=xt$, $y=yt$, $Y$ is given by \[
x^2t+y^2zt+2yz^m+\phi_{\geq 2}(yt,z,t)=0.\]
For $t=0$ we find that \[
f^{-1}(\Gamma)=E_1+mE_2 \]
and $E_1:y=t=0$, $E_2:z=t=0$. As before, $Y$ has exactly one singular point on $L=E_1 \cap E_2=(y,z,t)$. Let \[
\phi_{\geq 2}(y,z,t)=\phi_2 (y,z,t)+\phi_{\geq 3}(y,z,t).\]
Whether or not $W$ is terminal depends on what kind of singularities $Z$ has away from $E_2^{Z}$. In the chart $y=yt$, $Z$ 
is given by \[
x^2 +y^2t^2z+2yz^m+\frac{1}{t}\phi_2(yt^2,z,t)+\frac{1}{t}\phi_{\geq 3}(yt^2,z,t)=0.\]
Write again \[
\phi_{\geq 3}(yt^2,z,t)=tf_{\geq 2}(z)+t^2f_{\geq 1}(y,z,t) \]
with $f_{\geq 1}(y,0,0)=0$. Then $Z$ is given by \[
F=x^2+y^2t^2z+2yz^m+\frac{1}{t}\phi_2(yt^2,z,t)+f_{\geq 2}(z)+tf_{\geq 1}(y,z,t)=0.\]
Let $C=g^{-1}(0): x=z=t=0$, and \[
\Phi_2(y,z,t)=\frac{1}{t}\phi_2(yt^2,z,t).\]
We want to investigate the singularities of $Z$ along $C$. 
\begin{gather*}
\frac{\partial F}{\partial x}=2x  \\  
\frac{\partial F}{\partial y}=2yzt^2+2z^m+\frac{\partial \Phi_2}{\partial y}+t\frac{\partial f_{\geq 1}}{\partial y} \\
\frac{\partial F}{\partial z}=y^2t^2+2myz^{m-1}+\frac{\partial \Phi_2}{\partial z}+\frac{\partial f_{\geq 2}}{\partial z}
+t\frac{\partial f_{\geq 1}}{\partial z} \\
\frac{\partial F}{\partial t}=2y^2tz+\frac{\partial \Phi_2}{\partial t}+f_{\geq 1}(y,z,t)+t \frac{\partial f_{\geq 1}}{\partial t}
\end{gather*}
Hence $Z$ is singular along a point $Q\in C$ iff \[
\frac{\partial \Phi_2}{\partial y}(Q)=\frac{\partial \Phi_2}{\partial z}(Q)=\frac{\partial \Phi_2}{\partial t}(Q)=0. \]
Let \[
\phi_2(y,z,t)=a_1y^2+a_2t^2+a_3yt+a_4yz+a_5zt. \]
Then \[
\Phi_2(y,z,t)=\frac{1}{t}\phi_2(yt^2,z,t)=a_1y^2t^3+a_2t+a_3yt^2+a_4yzt+a_5z. \]
Now it follows that 
\begin{gather*}
\frac{\partial \Phi_2}{\partial y}=2a_1yt^3+a_3t^2+a_4zt \\
\frac{\partial \Phi_2}{\partial z}=a_4yt+a_5\\
\frac{\partial \Phi_2}{\partial t}=3a_1y^2t^2+a_2+2a_3yt+a_4yz
\end{gather*}
Along $C$, \[
\begin{array}{ccc}
\frac{\textstyle \partial \Phi_2}{\textstyle \partial y}=0,  & \frac{\textstyle \partial \Phi_2}{\textstyle \partial z}=a_5,   
& \frac{\textstyle \partial \Phi_2}{\textstyle \partial t}=a_2
\end{array} \]
Therefore, $Z$ is either singular along $C$ and $a_2=a_5=0$, or has exactly one singular point in $E_2^{Z}$ (in the other chart). 
If $a_2=a_5=0$, then \[
\phi_2(y,z,t)=y(a_1y+a_3t+a_4z).\]
A coordinate independent way to say this property is the statement of 2.b. 

We will now find the index of the singularities of $W$. Suppose that \[
K_{Z}=p^{\ast}K_W+aE_2^Z.\]
Moreover, $E_2^Z \cong E_2 \cong \mathbb{P}^2$. $E_1^{Z}$ is Cartier and therefore for a general line $l \subset E_2^Z$, $l \cdot E_1^Z = 1$. 
On the other hand, \[
l \cdot (E_1^Z+2E_2^Z)=l \cdot g^{\ast}(E_1+2E_2)=l \cdot (E_1 + 2E_2)=-1.\]
Hence $l \cdot E_2^Z=-1$. Moreover, \[
l \cdot K_Z = l \cdot g^{\ast}K_Y = l \cdot K_Y =-1.\]
Combining the above we see that $a=1$ and therefore \[
K_Z=p^{\ast}K_W+E_2^Z. \]
$E_2^Z$ has index 2 and hence $W$ has also index 2. From the above proof it is also clear that $W$ has exactly one index 2 point. This, 
as well as the type of the singularities, also follows from~\cite[Theorem 4.7]{Ko-Mo92}. 

Moreover, since $W$ has index 2 it follows 
that $-2E$ if $p-$very ample where $E=E_1^{W}$. The statement about the number of generators of $R(\Gamma, X)$ follows immediately.
\end{proof}
The difference between the $D_{2n}$ and $D_{2m+1}$ cases is shown by the next lemma.
\begin{lema}\label{diff}
Suppose that the general section $S$ of $X$ containing $\Gamma$ is $D_{2n+1}$ and $\Gamma \subset S$ is of type $DF_r$. 
Let $f: Y=B_{\Gamma}X \la X$. Let $f^{-1}(\Gamma)=E_1+dE_2$. Then $Y$ is singular along $L=E_1 \cap E_2$.
\end{lema}
This is precisely the reason that makes the $D_{2n+1}$ case very difficult to work with. 
In the $D_{2n}$ cases, $Y$ had exactly one singular point on $L$ and that made an explicit description of 
the $\mathbb{Q}$-factorialization of $E_1$ relatively easy.  
\begin{proof}
Suppose that $Y$ is not singular along $L$. Let $Q \in L$ be a smooth point. Let $S^{\prime}$ be a general section of $Y$ through $Q$. 
As in lemma 4.4, there is a section $S$ of $X$ through $\Gamma$ such that $S^{\prime}=f_{\ast}^{-1}S$. Then by assumption, 
$S$ is $D_{2n+1}$. Then by proposition 5.6.2.b.ii, it follows that $Q \in S^{\prime}$ is singular which is not true.
\end{proof}
The previous result could of course be proved by explicitely calculating $Y$ by using the normal forms of Proposition~\ref{normalforms}.

The results of Theorem~\ref{Dn} depend on the singularities of the general section $S$ of $X$ containing $\Gamma$. 
To apply the theorem it would be useful to get information about the general section from a special section. The next lemma 
gives informaion about the general section starting from a special one.
\begin{lema}
Let $\Gamma \subset X$. Suppose that the general section of $X$ containing $\Gamma$ is $D_k$. 
Let $P \in \Gamma \subset S_{0} \subset X$ be a special section and suppose that $S_{0}$ is $D_n$ with $n \geq 5$. Then
\begin{enumerate}
\item If $\Gamma \subset S_0$ is of type $DF_l$, then the general section $S$ of $X$ through $\Gamma$ 
is $D_m$, $m \leq n$ and also of type $DF_l$.
\item If $\Gamma \subset S_0$ is of type $DF_r$ and $n$ is even, then the 
general $S$ through $\Gamma$ is $D_{2k}$ and also of type $DF_r$.
\end{enumerate}
\end{lema}
\begin{proof}

\textbf{Case 1:} Suppose that $\Gamma \subset S_0$ is of type $FD_l$, i.e., 
$\Gamma^{\prime\prime}$ intersects $E_1$ in the fundamental cycle of $S_0$. Let $S$ be the general 
section through $\Gamma$ and assume it is $D_m$. If  $\Gamma^{\prime\prime}$ intersects $E_{m-1}$ or $E_m$ in the fundamental cycle of $S$, 
then by proposition 4.6 it follows that $d=m/2$, if $m$ is even, or $d=(m-1)/2$, if $m$ is odd. On the other hand, 
by the assumption on $S_0$ and Proposition~\ref{normalforms}, under suitable coordinates $\Gamma \subset X$ is given by \[
x^2+y^2z+z^{n-1}+t\phi_{\geq 2}(y,z,t)=0,\]
and $I_{\Gamma} =(x,z,t)$. Use notation as in lemma~\ref{diff}.

A computation as in Theorem~\ref{Dn} shows that $d=2$ and $Y$ has exactly one singular point on $L$. 
Hence the only possibility that the general section 
is not as claimed is that it is $D_5$ and $\Gamma^{\prime\prime}$ intersects $E_4$ or $E_5$. But then $Y$ is singular along $L$ 
as follows from lemma~\ref{diff}.

The fact that $m \leq n$ follows from the upper semicontinuity of the Tyurina number of the singularity.

\textbf{Case 2:} Suppose that $\Gamma \subset S_0$ is of type $FD_r$ and $n$ is even. 
First we will show that it is not possible that $\Gamma \subset S$, $S \subset X$ general through $\Gamma$, is of type $FD_l$. Suppose it is. 
Then since $n$ is even, $S^{\prime}=f_{\ast}^{-1}S$ has exactly one singular point which is $D_{n-1}$ as follows from proposition 4.6.2.a. 
On the other hand, $S_0^{\prime}$ has exactly one $A_{n-1}$ singular point $Q$. Therefore $Q\in Y$ is $cA_k$ and by~\cite{KoBa88} it 
is $cA_k$ in a neighborhood of $Q$. But then for a general $\Gamma \subset S \subset X$, $S^{\prime}$ is $A_k$ and hence it 
must be of type $FD_r$. 

If $\Gamma \subset S$ is $FD_r$ but $S$ is $D_{2m+1}$ for general $S$, then by lemma~\ref{diff} $Y$ is singular along $L$ which 
is not true as follows from corollary 4.7.
\end{proof}
Hence by looking at one section we know in which part of theorem~\ref{Dn} we are. 
So if we know that there is a section as in $1$ then all we need to know to conclude that there is no terminal 
contraction is that the general section is not $D_4$. The next lemma gives a criterion for that.
\begin{lema}
Let $\Gamma \subset X$ be given by \[
x^2+f_{\geq 3}(y,z,t)=0, \]
and $\Gamma =(x,y,t)$. Moreover suppose that $t=0$ is a DuVal section $S$ of $X$ containing $\Gamma$ and 
$\Gamma^{\prime\prime}$ intersects $E_1$ in the fundamental cycle of $S$. Then a $D_4$ section 
of $X$ containing $\Gamma$ does not exist iff \[
f_3(y,z,t) =  g(y,z,t) h^{2}(y,z,t).\]
\end{lema}
\begin{proof}
According to proposition~\ref{normalforms}, in suitable coordinates $\Gamma \subset X$ is given by \[
x^2+y^2z+z^{n-1}+t\phi_{\geq 2}(y,z,t)=0,\]
and $I_{\Gamma} =(x,z,t)$. The cubic term then of the above equation is \[
q_3(y,z,t)=y^2z+t\phi_2(y,z,t).\]
From the proof of the first part of proposition~\ref{normalforms} it follows that a $D_4$ section exists iff \[
q_3(y,z,t) \neq g(y,z,t)h^2(y,z,t), \]
for any $g(y,z,t)$, $h(y,z,t)$, which is the condition claimed by the lemma.
\end{proof}
\begin{example}
Let $X$ be given by \[
x^2+y^2z+2yz^n+t^m+t\phi_{\geq 3} (y,z,t)=0,\]
and $I_{\Gamma}=(x,y,t)$. Then there is a terminal contraction iff $m=3$. 
\end{example}
\begin{proof}
The section $S_0$ given by $(t=0)$ is $D_{2n}$ and $\Gamma \subset S_0$ is $DF_r$. Hence by lemma 6.3 
so is the general section through $\Gamma$, $S$. Now apply theorem~\ref{Dn}.
\end{proof}

Finally, I would like to mention that it will be interesting to get a more direct way of constructing the 
contraction of theorems~\ref{An},~\ref{Dn}. In particular, is it true that it is a weighted blow up of the curve? 
I plan to address this question as well as treat the remaining cases in a fture paper.